\documentclass[10pt,oneside]{jm}
\usepackage{amsfonts,amssymb,amsmath,url} 

\newcommand{\abs}{\vskip 0.5em\noindent}
\newcommand{\Abs}{\paragraph{}\hspace{-1em}}
\newcommand{\AbsT}[1]{\paragraph{\hspace{-1em} #1}}

\newcommand{\Pf}{\paragraph*{Proof.}}

\newcommand{\Def}{\AbsT{Definition.}}

\newcommand{\Rem}{\AbsT{Remark.}}

\newcommand{\Prop}{\AbsT{Proposition.}}

\newcommand{\bfi}{\noindent {\bf i)} }
\newcommand{\bfii}{\noindent {\bf ii)} }
\newcommand{\bfiii}{\noindent {\bf iii)} }

\newcommand{\bfa}{\noindent {\bf a)} }
\newcommand{\bfb}{\noindent {\bf b)} }
\newcommand{\bfc}{\noindent {\bf c)} }

\newcommand{\eg}{e.~g.\hskip 0.5em}
\newcommand{\ie}{i.~e.\hskip 0.5em}

\newcommand{\QED}{\hfill $\sharp$}

\newcommand{\N}{{\mathbb N}}
\newcommand{\Z}{{\mathbb Z}}
\newcommand{\Q}{{\mathbb Q}}

\newcommand{\C}{{\mathbb C}}
\newcommand{\F}{{\mathbb F}}

\newcommand{\cC}{{\mathcal C}}

\newcommand{\cE}{{\mathcal E}}

\newcommand{\cF}{{\mathcal F}}
\newcommand{\hcF}{{\hat{\mathcal F}}}
\newcommand{\cG}{{\mathcal G}}

\newcommand{\cM}{{\mathcal M}}

\newcommand{\cO}{{\mathcal O}}

\newcommand{\cR}{{\mathcal R}}
\newcommand{\cS}{{\mathcal S}}

\newcommand{\al}{\alpha}

\newcommand{\bt}{\beta}
\newcommand{\gm}{\gamma}

\newcommand{\lb}{\lambda}

\newcommand{\om}{\omega}

\newcommand{\Om}{\Omega}

\newcommand{\sg}{\sigma}

\newcommand{\chr}{\mbox{char}}

\newcommand{\End}{\mbox{End}}

\newcommand{\Hom}{\mbox{Hom}}

\newcommand{\id}{\mbox{id}}

\newcommand{\Stab}{\mbox{Stab}}

\newcommand{\ld}{,\ldots\hskip0em ,}
\newcommand{\lr}[1]{\langle #1\rangle}
\newcommand{\ov}[1]{\overline{#1}}

\newcommand{\wti}[1]{\widetilde{#1}}
\newcommand{\wh}[1]{\widehat{#1}}
\newcommand{\mt}{\mapsto}

\newcommand{\ra}{\rightarrow}

\newcommand{\cn}{\colon}

\newcommand{\sseq}{\subseteq}

\newcommand{\MAGMA}{{\sf MAGMA}}
\newcommand{\GAP}{{\sf GAP}}

\begin{document}
\raggedbottom
\setcounter{page}{1}
\pagestyle{myheadings}
\markboth{}{}
\thispagestyle{empty}
\begin{center} \Large\bf
On the ring of invariants of ordinary \\ 
quartic curves in characteristic $2$
\vspace*{1em} \\
\normalsize\rm
J\"urgen M\"uller and Christophe Ritzenthaler 
\footnote{C. Ritzenthaler acknowledges financial support provided
through the European Community's Human Potential Programme under
contract HPRN-CT-2000-00114, GTEM.} 
\vspace*{1em} \\
\end{center}


\section{Introduction}\label{intro}

\abs
Traditionally, non-singular curves of a fixed genus $g$ are classified 
in two categories, those which are hyperelliptic and those which are not. 
The first ones are usually considered to be the simplest, since their 
geometry relies on their Weierstraß points and is condensed on a line. 
Thus, the hyperelliptic locus $\cM_g^{h}$ of the corresponding moduli 
space $\cM_g$ is easily described and well understood. Seen from an
invariant theoretical viewpoint, this is reflected by the fact
that one only needs to deal with binary forms. To the contrary, 
non-hyperelliptic curves are more involved. Even in the simplest case,
non-singular quartic curves over $\C$, no complete description of the 
relevant invariant ring is known. More precisely, for the invariant ring 
of the natural action of $SL_3(\C)$ on the vector space of homogeneous 
polynomials of degree $4$ in $3$ variables a set of primary invariants,
see \cite{dixmier}, but no complete algebra generating set is known;
a conjecture of Shioda says that the invariant ring is generated
as an algebra by $13$ elements.

\abs
Hence the question arises whether we can describe the situation more  
precisely over other fields? If $F$ is a finite field of characteristic $2$, 
in \cite{nartritz} a complete classification of the $F$-isomorphism types 
of non-singular quartic curves defined over $F$ has been obtained.
Moreover, the stratification of the non-hyperelliptic locus $\cM_3^{nh}$
of the moduli space $\cM_3$ with respect to the $2$-rank of the Jacobian, 
and the $F$-rational points on the various strata, have been described there.
Here, the generic case is the one of ordinary non-singular quartic curves,
where the $2$-rank of the Jacobian is maximal, hence equal to $3$.
In \cite{nartritz}, a precise description of the ordinary non-singular 
quartic curves is given, and it is shown that the invariant ring associated 
to their moduli space $\cM_3^{ord}$ is given by a linear action of the 
finite group $G=GL_3(\F_2)$ on a $6$-dimensional vector space $W'^\ast$,
being defined over the field $\F_2$ with $2$ elements.
The aim of the present note is to give a complete and precise description
of this invariant ring, being denoted $S[W'^\ast]^G$ in the sequel,
thereby answering the corresponding question posed in \cite{nartritz}.

\abs
Despite the precise description of $G$ and $W'^\ast$, being suitable to be
handled by computer algebra systems,
brute force computer calculations to find primary and secondary 
invariants of $S[W'^\ast]^G$, using the standard press-button algorithms 
implemented in \MAGMA{} \cite{MAGMA}, as well as 
at most $2$ Gigabytes of memory and several hours of computing time, 
had to be abandoned unsuccessfully.
Hence the strategy employed here is to intertwine theoretical and 
computational analysis of $S[W'^\ast]^G$, which in effect leads both to some 
structural understanding of $S[W'^\ast]^G$ and finally to explicitly given 
invariants.
Actually, the theoretical analysis indicates how to combine ideas from 
computational invariant theory and tools already available in
computer algebra systems to obtain specially tailored techniques
applicable to the examples at hand.
Finally, the necessary computations have been carried out
the computer algebra systems \MAGMA{} \cite{MAGMA} and \GAP{} \cite{GAP}.
After all, to check the correctness of the results only needs a few 
seconds of computing time and approximately $10$ Megabytes of memory.
More details of the computations, a \MAGMA{} input file, as well as 
the primary and secondary invariants calculated, can be found on:
\url{http://www.math.jussieu.fr/~ritzenth}.

\abs
More precisely, the ring $S[W'^\ast]^G$ turns out to be Cohen-Macaulay, an 
optimal set of primary invariants has degrees $\{2,3,3,4,6,7\}$, and 
a corresponding minimal set of secondary invariants has cardinality $18$. 
Moreover, the algebra $S[W'^\ast]^G$ is generated by at most $12=6+6$ 
invariants, namely the $6$ primary invariants and $6$ of the secondary 
invariants, the latter having degrees $\{0,4,5,5,6,7\}$. Hence in 
particular $S[W'^\ast]^G$ is generated by invariants of degree at most $7$.
By the way, the number $12$ rings a bell: The authors wonder whether 
there is a connection to Shioda's conjecture mentioned above.

\abs
The paper is organized as follows:
In Section \ref{curves} we prepare the setting on ordinary quartic curves,
recall the necessary facts from \cite{nartritz}, and exhibit the $G$-module
$W'^\ast$ whose invariant ring $S[W'^\ast]^G$ we are going to examine.
In Section \ref{inv} we recall a few notions from commutative algebra and
derive some general facts about invariant rings needed in the sequel.
Finally, in Section \ref{ex} we carry out the analysis of the 
invariant ring $S[W'^\ast]^G$.
Note that we consider right group actions throughout, as this is common
is the computer algebra community and assumed in the computer algebra 
systems used.

\section{Ordinary quartic curves}\label{curves}

\Abs
Let $\F_2$ be the finite field with $2$ elements, and let
$\ov{\F}_2$ be its algebraic closure.
Let $M:=\F_2^3$ be the $3$-dimensional (row) 
vector space over $\F_2$, and let $\ov{M}:=M\otimes_{\F_2}\ov{\F}_2$.
Moreover, let $W:=\F_2^7$ and $\ov{W}:=W\otimes_{\F_2}\ov{\F}_2$, and let
$\cC:=\{C_{a,b,c,d,e,f,g};0\neq [a,b,c,d,e,f,g]\in\ov{W}\}$ 
be the family of quartic curves given by
$$ C_{a,b,c,d,e,f,g}\cn Q_{a,b,c,d,e,f}^2=g^2\cdot x y z(x+y+z),$$ 
where $Q_{a,b,c,d,e,f}:=ax^2+by^2+cz^2+dxy+eyz+fzx$.
This family is important because of the following

\Prop\label{sing} See \cite[Prop.1.1]{nartritz}.\\
Let $C$ be a non-singular quartic curve defined over $\ov{\F}_2$.
Then the following conditions are equivalent:

\bfi
The Jacobian variety $J_C$ of $C$ is ordinary, \ie we have
$|J_C[2](\ov{\F}_2)|=2^3$.

\bfii
The curve $C$ has seven bitangents.

\bfiii
The curve $C$ is isomorphic to some curve $C_{a,b,c,d,e,f,g}\in\cC$ 
such that

\abs\abs
$(\ast)$\hfill $abcg (a+b+d) (b+c+e) (a+c+f) (a+b+c+d+e+f+g)\neq 0$. \QED

\abs\abs\abs
Moreover, non-singular quartic curves $C$ and $C'$ are isomorphic,
if and only if there is an element $\gm\in GL_3(\ov{\F}_2)$ such that 
$C^\gm=C'$, where $\mathbb{P}(\ov{M})$ is considered as the natural
right module for the projective general linear group $PGL_3(\ov{\F}_2)$
of rank $3$ over $\ov{\F}_2$, and $C^\gm$ is the curve defined by
$$ V(C^\gm):=\{[x,y,z]\in\mathbb{P}(\ov{M});
               [x,y,z]\cdot\gm^{-1}\in V(C)\} ,$$
where $V(\cdot)$ denotes the locus of points of the curve.
Using the action on the $7$ bitangents, the isomorphism issue is 
reduced to the finite group $G:=GL_3(\F_2)=PGL_3(\F_2)\leq PGL_3(\ov{\F}_2)$. 
Indeed, non-singular curves 
$C_{a,b,c,d,e,f,g}\in\cC$ and $C_{a',b',c',d',e',f',g'}\in\cC$ 
are isomorphic, if and only if there is an element $\gm\in G$ such that 
$(C_{a,b,c,d,e,f,g})^\gm=C_{a',b',c',d',e',f',g'}$, see \cite{nartritz}.

\Abs\label{gdef}
In the sequel let $G:=GL_3(\F_2)$ be the general linear group 
of degree $3$ over $\F_2$, which up to isomorphism is the unique 
simple group of order $168$, see \cite[p.3]{Atlas}.
Let $A,B,C\in G$ be the elements of order $2$, $3$ and $7$, 
respectively, defined as 
$$ A:=\left[\begin{array}{ccc}
1 & . & 1 \\
. & 1 & . \\
. & . & 1 \\
\end{array}\right], \quad
B:=\left[\begin{array}{ccc}
. & . & 1 \\
1 & . & . \\
. & 1 & . \\
\end{array}\right], \quad
C:=\left[\begin{array}{ccc}
. & 1 & . \\
. & . & 1 \\
1 & 1 & . \\ \end{array}\right].$$

\abs
It is easily checked using \GAP{} or \MAGMA{}
that $G=\langle A,B\rangle$.
The $\F_2$-vector space $M$ can be 
considered as the natural right module for $G$. As the above
action of $PGL_3(\ov{\F}_2)$ restricts to an action of $G$ on $\cC$,
it is easily checked that we get an $\F_2$-linear action of $G$ on $W$ as
$$D_W\cn\quad A\mt\left[\begin{array}{ccc|ccc|c}
1 & . & . & . & . & . & . \\
. & 1 & . & . & . & . & . \\
1 & . & 1 & . & . & . & . \\
\hline
. & . & . & 1 & . & . & . \\
. & . & . & 1 & 1 & . & . \\
1 & . & . & . & . & 1 & . \\
\hline
. & . & . & 1 & . & . & 1 \\
\end{array}\right], \quad
B\mt\left[\begin{array}{ccc|ccc|c}
. & . & 1 & . & . & . & . \\
1 & . & . & . & . & . & . \\
. & 1 & . & . & . & . & . \\
\hline
. & . & . & . & . & 1 & . \\
. & . & . & 1 & . & . & . \\
. & . & . & . & 1 & . & . \\
\hline
. & . & . & . & . & . & 1 \\
\end{array}\right] .$$

\abs
The $\F_2$-subspace 
$W':=\{[a,b,c,d,e,f,g]\in W;g=0\}\leq W$ is a
$G$-submodule of $W$, and the above matrices show that
$W/W'\cong\F_2$ is the trivial $G$-module. Hence
we have an extension of $G$-modules

\abs\abs
$(\ast\ast)$ \hfill $\{0\}\ra W'\ra W\ra\F_2\ra\{0\}.$ \hfill \mbox{}

\abs\abs
Moreover, there is a $G$-submodule $W''\leq W'$, 
such that $\dim_{\F_2}(W'')=3$, as is also indicated above. 
The characteristic polynomials of the action of $C\in G$ on $W''$ and 
$W'/W''$ are $t^3+t+1\in\F_2[t]$ and $t^3+t^2+1\in\F_2[t]$,
respectively, which both are irreducible. Hence by \cite[p.3]{ModAtlas} 
we conclude that $W''$ and $W'/W''$ are non-isomorphic
absolutely irreducible $G$-modules, where $W''\cong M$ is isomorphic to 
the natural representation of $G$, and $W'/W''$
is obtained from $W''$ by applying the automorphism of $G$ given
by inverting and transposing matrices.

\Abs\label{ext}
Let $\eta\in Z^1(G,W')$ be the cocycle describing the extension $(\ast\ast)$,
and let $\cS_4\cong H\leq G$ be the subgroup permuting the set
$\{x^\ast,y^\ast,z^\ast,(x+y+z)^\ast\}\sseq M^\ast$, where 
$M^\ast:=\Hom_{\F_2}(M,\F_2)$ is the $G$-module contragredient to $M$
and $\{x^\ast,y^\ast,z^\ast\}\sseq M^\ast$ is
the $\F_2$-basis dual to the standard basis of $M$.

\abs
By construction, for the restriction of $\eta$ to the subgroup
$H$ we have $\eta|_H=0\in Z^1(H,W')$.
As $[G\cn H]=7$ is invertible in $\F_2$,
we by \cite[Cor.3.6.18]{benson} conclude that
$\eta=0\in H^1(G,W')\cong\mbox{Ext}_G^1(\F_2,W')$. Thus the
extension $(\ast\ast)$ splits, and we have $W=W'\oplus\F_2$ as $G$-modules.
More concretely, going over from the standard basis of $W$ 
to the $\F_2$-basis where the last standard basis vector is replaced by
$[1,1,1,1,1,1,1]\in W$, we indeed obtain
$$D'_W\cn\quad A\mt\left[\begin{array}{ccc|ccc|c}
1 & . & . & . & . & . & . \\
. & 1 & . & . & . & . & . \\
1 & . & 1 & . & . & . & . \\
\hline
. & . & . & 1 & . & . & . \\
. & . & . & 1 & 1 & . & . \\
1 & . & . & . & . & 1 & . \\
\hline
. & . & . & . & . & . & 1 \\
\end{array}\right], \quad
B\mt\left[\begin{array}{ccc|ccc|c}
. & . & 1 & . & . & . & . \\
1 & . & . & . & . & . & . \\
. & 1 & . & . & . & . & . \\
\hline
. & . & . & . & . & 1 & . \\
. & . & . & 1 & . & . & . \\
. & . & . & . & 1 & . & . \\
\hline
. & . & . & . & . & . & 1 \\
\end{array}\right] .$$

\abs
Actually, this basis change amounts to 
substituting the curves $C_{a,b,c,d,e,f,g}\in\cC$ by curves
defined by
$$ (a x^2 +b y^2 +c z^2+d xy+eyz+f zx)^2=g^2\cdot C_K(x,y,z) ,$$
where $C_K(x,y,z)=x^4+y^4+z^4+(xy)^2+(yz)^2+(zx)^2+xyz(x+y+z)$. 
We have 
$x^{\ast 4}+y^{\ast 4}+z^{\ast 4}
 +(x^\ast y^\ast)^2+(y^\ast z^\ast)^2+(z^\ast x^\ast)^2
 +x^\ast y^\ast z^\ast (x^\ast+y^\ast+z^\ast)\in S[M^\ast]^G$
by construction, where $S[M^\ast]^G\sseq S[M^\ast]$ 
denotes the $\F_2$-subalgebra of $G$-invariants of the symmetric 
algebra $S[M^\ast]$ over $M^\ast$. 
As for the homogeneous component $S[M^\ast]^G_4$ of $S[M^\ast]^G$ 
of degree $4$ we have $\dim_{\F_2}(S[M^\ast]^G_4)=1$, see \cite{elkies}, 
the curve $C_K\sseq\mathbb{P}(\ov{M})$ is a twist of the Klein quartic curve.

\abs
As the extension $(\ast\ast)$ splits, for the corresponding contragredient
$G$-modules we have $W^\ast\cong W'^\ast\oplus\F_2^\ast$. 
With respect to the $\F_2$-basis
$\{a^\ast,b^\ast,c^\ast,d^\ast,e^\ast,f^\ast\}\sseq W'^\ast$
dual to the standard basis of $W'$, the $G$-action on $W'^\ast$ is given as
$$D_{W'^\ast}\cn\quad A\mt\left[\begin{array}{ccc|ccc}
1 & . & 1 & . & . & 1 \\
. & 1 & . & . & . & . \\ 
. & . & 1 & . & . & . \\
\hline
. & . & . & 1 & 1 & . \\ 
. & . & . & . & 1 & . \\
. & . & . & . & . & 1 \\
\end{array}\right], \quad
B\mt\left[\begin{array}{ccc|ccc}
. & . & 1 & . & . & . \\
1 & . & . & . & . & . \\
. & 1 & . & . & . & . \\
\hline
. & . & . & . & . & 1 \\
. & . & . & 1 & . & . \\
. & . & . & . & 1 & . \\
\end{array}\right] .$$

\abs
Let $\ov{W'^\ast}:=W'^\ast\otimes_{\F_2}\ov{\F}_2$ and 
$\ov{W'}:=W'\otimes_{\F_2}\ov{\F}_2$. Moreover, let
$S[W'^\ast]$ denote the symmetric $\F_2$-algebra over $W'^\ast$.
Hence we have $S[\ov{W'^\ast}]\cong S[W'^\ast]\otimes_{\F_2}\ov{\F}_2$.
Moreover, by Remark \ref{fieldextrem} below we have
$S[\ov{W'^\ast}]^G\cong S[W'^\ast]^G\otimes_{\F_2}\ov{\F}_2$.
The embedding of affine rings $S[\ov{W'^\ast}]^G\sseq S[\ov{W'^\ast}]$, 
where $S[\ov{W'^\ast}]^G$ is the ring of $G$-invariants in $S[\ov{W'^\ast}]$,
defines a morphism $\ov{W'}\ra\ov{W'}/G$ of affine varieties 
over $\ov{\F}_2$, which as $G$ is finite is a geometric quotient,
see \cite[Ch.2.3]{derksen}. Hence we have proved the following

\Prop See \cite[Prop.1.3]{nartritz}. \\
The moduli space $\cM_3^{ord}$ of the ordinary quartic curves
is isomorphic to the open subset of the affine variety 
$\mbox{Spec}(S[W'^\ast]^G\otimes_{\F_2}\ov{\F}_2)$
given by the non-singularity conditions $(\ast)$, 
see Proposition \ref{sing}.
\QED

\section{Invariant rings}\label{inv}

\abs
We recall a few notions from commutative algebra and
derive some general facts about invariant rings;
as general references see \eg \cite{benson2,derksen}.
Note that we use the following piece of notation frequently:
If $\cE$ is a subset of a commutative ring $R$, then 
$\cE R$ denotes the ideal of $R$ 
generated by $\cE$.

\Def
Let $R=\bigoplus_{d\geq 0}R_d$ be a finitely generated commutative 
$\Z^{\geq 0}$-graded algebra over a field $F$, such that 
$\dim_F(R_d)<\infty$ for $d\in\N_0$, and $R_0\cong F$. 
Let $R_+:=\bigoplus_{d>0}R_d\lhd R$ be the {\bf irrelevant ideal},
and let  $H_R(t)\in\Q(t)\sseq \Q((t))$ denote the 
{\bf Hilbert series} of $R$.

\abs
A set $\cF:=\{f_1\ld f_n\}\sseq R_+$ of homogeneous elements,
where $n=\dim(R)$ is the {\bf Krull dimension} of $R$,
is called a {\bf homogeneous system of parameters}, if $\cF$ is 
algebraically independent and if $R$ is a finitely generated 
$F[\cF]$-module, where $F[\cF]\sseq R$ is the polynomial 
$F$-subalgebra of $R$ generated by $\cF$.

\Def
The $F$-algebra $R$ is called {\bf Cohen-Macaulay} if there is a 
homogeneous system of parameters $\cF$ such that $R$ is a free 
$F[\cF]$-module.
If $R$ is Cohen-Macaulay, then $R$ is a free $F[\cF]$-module 
for each homogeneous system of parameters $\cF$. 
Moreover, $R$ is Cohen-Macaulay if and only if $\dim(R)=\mbox{depth}(R)$.
Here the {\bf depth} of $R$ is the common length of 
all maximal {\bf regular} homogeneous sequences in $R$, 
see \cite[Ch.4.3]{benson2}.
Note that each regular sequence can be extended to a maximal one; and
if $R$ is Cohen-Macaulay, then the homogeneous systems of parameters of 
$R$ and the maximal regular homogeneous sequences in $R$ coincide.

\Rem\label{hilserrem}
Let $R$ be Cohen-Macaulay, and let $\cF=\{f_1\ld f_n\}\sseq R$ 
be a homogeneous system of parameters.
Hence $R$ is a free $F[\cF]$-module, and 
for a minimal homogeneous $\cF$-module generating set 
$\cG=\{g_1\ld g_s\}$ of $R$ we have
$$ H_R(t)=\frac{\sum_{j=1}^s t^{e_j}}{\prod_{i=1}^n (1-t^{d_i})}\in\Q(t) ,$$
where $d_i=\deg(f_i)$ and $e_j= \deg(g_j)$.
Hence in the Cohen-Macaulay case the degrees $d_i$ and $e_j$ 
of the elements of $\cF$ and $\cG$ can be read off from $H_R$.

\Def
Let $G$ be a finite group, let $V$ be an $FG$-module
and let $S[V]^G$ denote the ring of {\bf $G$-invariants} in
the symmetric algebra $S[V]$ over $V$.
By \cite[Thm.1.3.1]{benson2} 
the $F$-algebra $S[V]^G$ is finitely generated and we
have $\dim(S[V]^G)=\dim(S[V])=\dim_F(V)$.
A homogeneous system of parameters $\cF$ of $S[V]^G$ is called a set of 
{\bf primary invariants}. A minimal set of homogeneous $F[\cF]$-module 
generators of $R$ is called a set of {\bf secondary invariants}.

\Rem\label{fieldextrem}
Let $F\sseq L$ be a field extension, and let $V_L:=V\otimes_F L$.
As for $d\in\N_0$ we have
$S[V]^G_d=\bigcap_{\sg\in G}\ker_{S[V]_d}(\sg -1)$, 
we conclude that $\dim_F(S[V]^G_d)=\dim_L(S[V_L]^G_d)$, for $d\in\N_0$. 
Thus we have $S[V]^G\otimes_F L=S[V_L]^G$.

\Def
Let $H\leq G$ be a subgroup such that the characteristic $\chr(F)$ 
of the field $F$ does not divide the index $[G\cn H]$ of $H$ in $G$.
Then the {\bf relative Reynolds operator} with respect to $H$ and $G$
is defined as
$$ \cR_H^G\cn S[V]^H \ra S[V]^G\cn 
              f \mt\frac{1}{[G\cn H]}\cdot\sum_{\sg\in H|G}f\cdot\sg ,$$
where the sum runs over a set of representatives of the right cosets
$H|G$ of $H$ in $G$. 
Note that we have $S[V]^G\sseq S[V]^H$, and hence
$\cR_H^G$ is an $S[V]^G$-module projection onto $S[V]^G$.

\abs\abs
The following Proposition is a slight generalization of the 
Hochster-Eagon Theorem, see \cite[Thm.4.3.6]{benson2}, saying that
if $\chr(F)$ does not divide $|G|$, 
then $S[V]^G$ is Cohen-Macaulay; see also \cite[Prop.8.3.1]{smith}.

\Prop\label{cmsubgprop}
Let $H\leq G$ be a subgroup such that $\chr(F)$ does not divide $[G\cn H]$.
Provided $S[V]^H$ is Cohen-Macaulay, then $S[V]^G$ is Cohen-Macaulay as well.

\Pf
Let $\cF\sseq S[V]^G$ be a set of primary invariants of $S[V]^G$.
As both $F[\cF]\sseq S[V]^G$ and $S[V]^G\sseq S[V]^H$ are
finite ring extensions, this also holds for $F[\cF]\sseq S[V]^H$,
and hence $\cF\sseq S[V]^H$ is a set of primary invariants of $S[V]^H$.
As $S[V]^H$ is Cohen-Macaulay, it hence is a free $F[\cF]$-module.
Using the relative Reynolds operator $\cR_H^G$, we conclude that
$S[V]^G$ is a direct summand of the graded $F[\cF]$-module $S[V]^H$.
Hence $S[V]^G$ is a finitely generated graded projective $F[\cF]$-module, 
and thus by \cite[La.4.1.1]{benson2} is a free $F[\cF]$-module.
\QED

\Rem\label{primrem}
To compute primary invariants in our particular situation in
Section \ref{ex}, we will exploit the following setting.

\abs
Let $\{0\}\ra U'\stackrel{\al}{\ra} U\stackrel{\bt}{\ra}V\ra\{0\}$ 
be an extension of $G$-modules. Hence $\al$ induces an embedding 
$S[U']\sseq S[U]$, and $\bt$ induces an isomorphism 
$S[U]/U'S[U]\ra S[V]$. 
Hence we have $S[U']^G\sseq S[U]^G$ and 
$(S[U]/U'S[U])^G\cong S[V]^G$, 
where $U'S[U]\lhd S[U]$ is a $G$-submodule. Unfortunately,
in general we only have an embedding 
$S[U]^G/(U'S[U])^G\sseq (S[U]/U'S[U])^G$, but not an 
isomorphism, and in general $(U'S[U])^G\lhd S[U]^G$ 
is not generated by $S[U']_+^G$.

\abs
Let us assume that the above extension splits, and let 
$\gm\cn V\ra U$ such that $\gm\bt=\id_{V}$.
Hence $\gm$ induces an embedding $S[V]\sseq S[U]$, and we have 
$S[U]=S[V]\oplus U'S[U]$ as $G$-modules. Thus from
$S[U]^G=S[V]^G\oplus (U'S[U])^G$ we conclude that $\bt$ 
induces an isomorphism $S[U]^G/(U'S[U])^G\cong S[V]^G$.

\abs
Let us moreover assume that $U'\cong F$ is the trivial $G$-module,
and let $0\neq \hat{f}\in U'\sseq S[U]^G_1$.
Then we have $(U'S[U])^G=\hat{f}\cdot S[U]^G\lhd S[U]^G$, 
and thus $\bt$ induces an isomorphism $S[U]^G/\hat{f}S[U]^G\cong S[V]^G$.
As $\hat{f}\in S[U]$ is not a zero-divisor, 
for the corresponding Hilbert series 
we obtain $H_{S[V]^G}(t)=(1-t)\cdot H_{S[U]^G}(t)\in\Q(t)$.

\abs
Let $\hcF:=\{\hat{f}_0\ld\hat{f}_{n-1}\}\sseq S[U]^G$ 
be a set of primary invariants, such that $\hat{f}_0=\hat{f}$,
and let $f_i:=\hat{f}_i\bt\in S[V]^G$, for $i\in\{1\ld n-1\}$,
as well as $\cF:=\{f_1\ld f_{n-1}\}$.
Note that we have $n=\dim(S[U]^G)=\dim_F(U)$ and
$n-1=\dim(S[V]^G)=\dim_F(V)$.
As $\hcF\sseq S[U]^G$ is a set of primary invariants,
by the Graded Nakayama Lemma, see \cite[La.3.5.1]{derksen}, 
we have $\dim(S[U]^G/\hcF S[U]^G)=0$.
Hence using $\bt$ we find $\dim(S[V]^G/\cF S[V]^G)=0$,
where $|\cF|\leq n-1$. Hence by the Graded Nakayama Lemma again
we conclude that $\cF\sseq S[V]^G$ is set of primary invariants.

\abs
Finally, we note the following. As $\hat{f}\in S[U]^G$ is not a 
zero-divisor, and hence regular, we conclude that $S[U]^G$ is 
Cohen-Macaulay, if and only if $S[V]^G$ is. Moreover,
if $S[V]^G$ is Cohen-Macaulay, then there is a maximal regular
homogeneous sequence in $S[U]^G$ beginning with $\hat{f}$, which
thus is a set of primary invariants of $S[U]^G$, and the above 
construction indeed yields a a set of primary invariants of $S[V]^G$,
which is optimal in the sense of \cite[Ch.3.3.2]{derksen} if the used
set of primary invariants of $S[U]^G$ was.

\Rem\label{secrem}
To compute secondary invariants in our particular situation 
in Section \ref{ex}, we use a special adaptation of the method 
typically used in the case where $\chr(F)$ does not divide $|G|$, 
see \cite[Ch.3.5]{derksen}.

\abs
Let us assume that $R:=S[V]^G$ is known to be Cohen-Macaulay,
and that the Hilbert series $H_R\in\Q(t)$ and a set $\cF\sseq R_+$
of primary invariants are known.
Hence by Remark \ref{hilserrem} we have 
$f:=\prod_{i=1}^n(1-t^{d_i})\cdot H_R\in\Z^{\geq 0}[t]$,
and hence the cardinality $s$ of any minimal homogeneous $F[\cF]$-module 
generating set $\cG$ of $R$ is given as $s=f(1)$, while the degrees $e_j$ 
can be determined from the monomials occurring in $f$.

\abs
Let $\cG\sseq R$ be a set having the appropriate cardinality, 
whose elements are homogeneous of the appropriate degrees.
By the Graded Nakayama Lemma, see \cite[La.3.5.1]{derksen},
the set $\cG$ generates the $F[\cF]$-module $R$, if and only if 
$\cG$ generates the $F$-vector space $R/F[\cF]_+ R$.
By the assumptions made on $\cG$ we conclude that 
$\cG$ is a generating set of the $F[\cF]$-module $R$,
if and only if $\cG\sseq R/F[\cF]_+ R$ is $F$-linearly independent.

\abs
As we are developing a method to find secondary invariants, the ring $R$ and
hence $R/F[\cF]_+ R$ are not yet known. Thus  we proceed as follows.
Let $H\leq G$ be a subgroup such that $\chr(F)$ does not divide 
$[G\cn H]$, and let $S:=S[V]^H$. As we have $F[\cF]\sseq R\sseq S$,
we may consider the natural map $\pi\cn R\ra S\ra S/(\sum_{i=1}^n f_iS)$
of $F$-algebras. Hence we have $F[\cF]_+ R\sseq\ker(\pi)$. 
Conversely, let $h\in\ker(\pi)\unlhd R$, hence we have
$h=\sum_{i=1}^n f_i h_i$, where $h_i\in S$. Thus we have
$h=\cR_H^G(h)=\sum_{i=1}^n f_i\cdot\cR_H^G(h_i)\in F[\cF]_+ R$, and
hence $\ker(\pi)=F[\cF]_+ R$. Thus we have an embedding
$\pi\cn R/F[\cF]_+ R\ra S/(\sum_{i=1}^n f_iS)$. Hence
$\cG\sseq R/F[\cF]_+ R$ is $F$-linearly independent,
if and only if $\pi(\cG)\sseq S/(\sum_{i=1}^n f_iS)$ is.

\abs
Let us finally assume that $S=S[V]^H$ is Cohen-Macaulay, and that
$\cG'\sseq S$ is a minimal set of secondary invariants; note
that $\cF\sseq S$ is a set of primary invariants.
As $S$ is the free $F[\cF]$-module generated by $\cG'$, a description
of $S$ as a finitely presented commutative $F$-algebra can be derived
using linear algebra techniques, see \cite[Ch.3.6]{derksen}.
From that, a description of $S/(\sum_{i=1}^n f_iS)$ as a finitely
presented commutative $F$-algebra is immediately derived. Hence 
the $F$-linear independence of $\pi(\cG)\sseq S/(\sum_{i=1}^n f_iS)$ 
is easily verified or falsified using Gröbner basis techniques,
see \cite[Ch.3.5]{derksen}.
\QED

\section{The invariant ring $S[W'^\ast]^G$}\label{ex}

\abs
Let again $G:=GL_3(\F_2)$ and let $S[W'^\ast]^G$ be
the invariant ring introduced in Section \ref{ext}. 
We are prepared to analyze its structure, and begin with a
module theoretic property of the $G$-module $W$.

\Prop\label{tsmodprop}
The $G$-module $W$ is a transitive permutation module.

\Pf
If $W$ were a transitive permutation module, 
the corresponding point stabilizer would be a subgroup of order $24$, 
leading to the following sensible guess.
As in Section \ref{ext}, let $\cS_4\cong H\leq G$ be the subgroup 
permuting the set $\{x^\ast,y^\ast,z^\ast,(x+y+z)^\ast\}\sseq M^\ast$.
Hence $H$ fixes $w_0=[0,0,0,0,0,0,1]\in W$, and we are led to conjecture 
that $H=\Stab_G(w_0)$ and that $w_0\cdot G\sseq W$ is an $\F_2$-basis 
of $W$ being permuted by $G$.
To check this, let $C\in G$ be as in Section \ref{gdef}.
Its action on $W$ is given as
$$D_W\cn\quad C\mt\left[\begin{array}{ccc|ccc|c}
. & . & 1 & . & . & . & . \\
1 & . & 1 & . & . & . & . \\
. & 1 & . & . & . & . & . \\
\hline
. & . & 1 & . & . & 1 & . \\
. & . & . & 1 & 1 & . & . \\
. & . & . & . & 1 & . & . \\
\hline
. & . & . & . & 1 & . & 1 \\
\end{array}\right] .$$

\abs
Hence $\{C^i\in G;i=0,\ldots,6\}\sseq G$ is a set of representatives 
of the right cosets $H|G$, and 
$\Om:=\{w_0\cdot C^i\in W;i=0,6,1,2,3,4,5\}\sseq W$ 
is given as follows, where the rows indicate the elements of $\Om$ 
in terms of the standard basis of $W$,
$$\left[\begin{array}{ccccccc}
 . & . & . & . & . & . & 1 \\
 . & . & . & . & . & 1 & 1 \\ 
 . & . & . & . & 1 & . & 1 \\
 . & . & . & 1 & . & . & 1 \\
 . & . & 1 & . & 1 & 1 & 1 \\
 . & 1 & . & 1 & 1 & . & 1 \\
 1 & . & . & 1 & . & 1 & 1 \\
\end{array}\right] .$$

\abs
Hence $\Om$ is an $\F_2$-basis of $W$, and it is easily checked 
that it is permuted by $G$, where in particular
$A\mt (1,4)(2,7)\in\cS_7$ and $B\mt (2,4,3)(5,7,6)\in\cS_7$.
\QED

\abs\abs
By Proposition \ref{tsmodprop} the $G$-module $W'^\ast$ is a 
direct summand of the permutation module $W^\ast$.
Note that direct summands of permutation modules are also called 
{\bf trivial source modules}, see \cite[Ch.3.10.4]{derksen}.

\Prop\label{hilprop}
The Hilbert series $H_{S[W'^\ast]^G}\in\Q(t)$ of $S[W'^\ast]^G$ is given as
$$ H_{S[W'^\ast]^G}(t)=
\frac{1\!+\!t^4\!+\!2 t^5\!+\!t^6\!+\!t^7\!+\!t^8\!+\!2 t^9
      \!+\!2 t^{10}\!+\!t^{11}\!+\!t^{12}\!+\!t^{13}\!+\!
      2 t^{14}\!+\!t^{15}\!+\!t^{19}}
     {(1-t^2)\cdot (1-t^3)^2\cdot (1-t^4)\cdot (1-t^6)\cdot (1-t^7)}.$$

\Pf
As the $G$-module $W'^\ast$ is a trivial source module,
by \cite[Cor.3.11.4]{benson} the $G$-module $W'^\ast$ has a 
unique lift to a trivial source $\Z_2 G$-module $\wh{W'^\ast}$, 
where $\Z_2\sseq \Q_2$ is the integral closure of $\Z$ in the
$2$-adic completion $\Q_2$ of $\Q$. 
Let $\wti{W'^\ast}:=\wh{W'^\ast}\otimes_{\Z_2}\Q_2$.
By \cite[Prop.3.10.15]{derksen} we have
$H_{S[W'^\ast]^G}(t)=H_{S[\wti{W'^\ast}]^G}(t)\in\Q(t)$,
where by Molien's Theorem, see \cite[Thm.2.5.2]{benson2},
the latter is given as
$$ H_{S[\wti{W'^\ast}]^G}(t)=
  \frac{1}{|G|}\cdot\sum_{\sg\in G} 
  \frac{1}{\det_{\wti{W'^\ast}}(1-t\sg)} \in\Q_2(t) .$$
Moreover, we have 
$\det_{\wti{W'^\ast}}(1-t\sg)=\prod_{i=1}^7 (1-\lb_i(\sg)\cdot t)$,
where 
$\{\lb_1(\sg),\ldots,\lb_7(\sg)\}$
are the eigenvalues of the action of $\sg\in G$ in a suitable extension
field of $\Q_2$. Hence $\det_{\wti{W'^\ast}}(1-t\sg)$ can be evaluated
from the ordinary character table of $G$, see \cite[p.3]{Atlas}, 
and the character $\chi_{\wti{W'^\ast}}$ of $\wti{W'^\ast}$, 
see \cite[Ch.2.5]{benson2}. This method is implemented in \GAP,
where also the ordinary character table of $G$ is available.

\abs
Hence it remains to find the character $\chi_{\wti{W'^\ast}}$.
Thus we have to determine the trivial source lift $\wh{W'^\ast}$. 
By the proof of Proposition \ref{tsmodprop} we
have $(\F_2)_H^G\cong W\cong W'\oplus\F_2$, where $\F_2$ denotes the 
trivial $G$-module. Hence we have the trivial source lifts
$(\Z_2)_H^G\cong\wh{W}\cong\wh{W'}\oplus\Z_2$, where 
again $\Z_2$ denotes the trivial $G$-module. 
Since we have $(\wh{W})^\ast\cong\wh{W^\ast}$ as $G$-modules,
we obtain $((\Z_2)_H^G)^\ast\cong\wh{W^\ast}\cong\wh{W'^\ast}\oplus\Z_2$,
Moreover, we conclude
$((\Q_2)_H^G)^\ast\cong\wh{W^\ast}\otimes_{\Z_2}\Q_2
 \cong\wti{W'^\ast}\oplus\Q_2$.
By \cite[p.3]{Atlas} the character $1_H^G$ of the permutation 
$G$-module $(\Q_2)_H^G$ is given as $1_G+\chi_6$, where $\chi_6$ is
the unique irreducible character of degree $6$ and $1_G$ is
the trivial character. As $\chi_6$ is real-valued, we hence have 
$\chi_{\wti{W'^\ast}}=1_H^G-1_G=\chi_6$.
\QED

\Prop\label{cmprop}
The invariant ring $S[W'^\ast]^G$ is Cohen-Macaulay.

\Pf
Let $D\leq G$ be a $2$-Sylow subgroup of $G$, hence we have $|D|=8$.
Using the standard methods to compute primary invariants, 
see \cite[Ch.3.3]{derksen}, and secondary invariants in the modular case,
see \cite[Ch.3.5]{derksen}, implemented in \MAGMA, 
we find primary invariants $\{f'_1\ld f'_6\}\sseq S[W'^\ast]^D$ 
having degrees $\{1,1,2,2,2,4\}$, and a minimal set of secondary invariants
$\cG':=\{g'_0\ld g'_3\}\sseq S[W'^\ast]^D$ having degrees 
$\{0,3,3,6\}$, where of course $g'_0=1$. 
As we moreover have $|\cG'|\cdot |D|=\prod_{i=1}^6 \deg(f'_i)$, 
by \cite[Thm.3.7.1]{derksen} we conclude that $S[W'^\ast]^D$ 
is Cohen-Macaulay, and thus by Proposition \ref{cmsubgprop} 
the ring $S[W'^\ast]^G$ also is.
\QED

\Abs
We are prepared to compute primary invariants of $S[W'^\ast]^G$.
To do this, we first consider 
the permutation module $W^\ast=W'^\ast\oplus\F_2$, 
and compute the homogeneous components $S[W^\ast]^G_d$, for $d\leq 7$, 
as follows. Let $\Om^\ast=\{\om_1^\ast\ld \om_7^\ast\}\sseq W^\ast$ 
be the $\F_2$-basis of $W^\ast$ dual to the $\F_2$-basis 
$\Om\sseq W$ given in the proof of Proposition \ref{tsmodprop}.
As $S[W^\ast]_d$ also is a permutation module, whose $\F_2$-basis 
$(\Om^\ast)^d$, consisting of the monomials of degree $d$
in the indeterminates $\Om^\ast$, is permuted by $G$.
Hence $(\Om^\ast)^d$ is partitioned into $G$-orbits
$(\Om^\ast)^d=\coprod_{i=1}^{n_d}\cO_i$, 
where $n_d=\dim_{\F_2}(S[W^\ast]^G_d)$. Letting 
$\cO_i^+:=\sum_{f\in\cO_i}f\in S[W^\ast]^G_d$
denote the corresponding orbit sum, the set
$\{\cO_i^+;i\in\{1\ld n_d\}\}$ forms an 
$\F_2$-basis of $S[W^\ast]^G_d$, see also \cite[Ch.3.10]{derksen}.

\abs
As by Remark \ref{primrem} we have 
$H_{S[W^\ast]^G}=\frac{1}{1-t}\cdot H_{S[W'^\ast]^G}\in\Q(t)$,
we look for primary invariants having degrees $\{1,2,3,3,4,6,7\}$.
By \cite[Prop.3.3.1]{derksen}, a set 
$\hcF=\{\hat{f}_0\ld\hat{f}_6\}\sseq S[W^\ast]^G$ 
of homogeneous elements is a set of primary invariants, if and only if
$\dim(S[W^\ast]/\hcF S[W^\ast])=0$.
Krull dimensions can be computed using Gröbner basis techniques,
which are implemented in \MAGMA{}, and we indeed find the 
following set $\hcF$ of primary invariants of $S[W^\ast]^G$, 
consisting of certain of the orbit sums computed above,

$$ \begin{array}{rcl}
\hat{f}_0 & := & 
\om^\ast_1+\om^\ast_2+\om^\ast_3+\om^\ast_4
+\om^\ast_5+\om^\ast_6+\om^\ast_7, \\ 
\hat{f}_1 & := & 
\om^\ast_1\om^\ast_2+\om^\ast_1\om^\ast_3+\om^\ast_1\om^\ast_4
+\om^\ast_1\om^\ast_5+\om^\ast_1\om^\ast_6+\om^\ast_1\om^\ast_7
+\om^\ast_2\om^\ast_3+ \\
& &
\om^\ast_2\om^\ast_4+\om^\ast_2\om^\ast_5+\om^\ast_2\om^\ast_6
+\om^\ast_2\om^\ast_7+\om^\ast_3\om^\ast_4+\om^\ast_3\om^\ast_5
+\om^\ast_3\om^\ast_6+ \\
& &
\om^\ast_3\om^\ast_7+\om^\ast_4\om^\ast_5+\om^\ast_4\om^\ast_6
+\om^\ast_4\om^\ast_7+\om^\ast_5\om^\ast_6+\om^\ast_5\om^\ast_7
+\om^\ast_6\om^\ast_7, \\ 
\hat{f}_2 & := & 
\om^\ast_1\om^\ast_2\om^\ast_3+\om^\ast_1\om^\ast_2\om^\ast_4
+\om^\ast_1\om^\ast_2\om^\ast_5+\om^\ast_1\om^\ast_2\om^\ast_7
+\om^\ast_1\om^\ast_3\om^\ast_4+\om^\ast_1\om^\ast_3\om^\ast_5+\\
& &
\om^\ast_1\om^\ast_3\om^\ast_6+\om^\ast_1\om^\ast_4\om^\ast_6
+\om^\ast_1\om^\ast_4\om^\ast_7+\om^\ast_1\om^\ast_5\om^\ast_6
+\om^\ast_1\om^\ast_5\om^\ast_7+\om^\ast_1\om^\ast_6\om^\ast_7+\\
& &
\om^\ast_2\om^\ast_3\om^\ast_5+\om^\ast_2\om^\ast_3\om^\ast_6
+\om^\ast_2\om^\ast_3\om^\ast_7+\om^\ast_2\om^\ast_4\om^\ast_5
+\om^\ast_2\om^\ast_4\om^\ast_6+\om^\ast_2\om^\ast_4\om^\ast_7+\\
& &
\om^\ast_2\om^\ast_5\om^\ast_6+\om^\ast_2\om^\ast_6\om^\ast_7
+\om^\ast_3\om^\ast_4\om^\ast_5+\om^\ast_3\om^\ast_4\om^\ast_6
+\om^\ast_3\om^\ast_4\om^\ast_7+\om^\ast_3\om^\ast_5\om^\ast_7+\\
& &
\om^\ast_3\om^\ast_6\om^\ast_7+\om^\ast_4\om^\ast_5\om^\ast_6
+\om^\ast_4\om^\ast_5\om^\ast_7+\om^\ast_5\om^\ast_6\om^\ast_7, \\
\hat{f}_3 & := &
\om^\ast_1\om^\ast_2\om^\ast_6+\om^\ast_1\om^\ast_3\om^\ast_7
+\om^\ast_1\om^\ast_4\om^\ast_5+\om^\ast_2\om^\ast_3\om^\ast_4+\\
& &
\om^\ast_2\om^\ast_5\om^\ast_7+\om^\ast_3\om^\ast_5\om^\ast_6
+\om^\ast_4\om^\ast_6\om^\ast_7, \\
\hat{f}_4 & := &
\om^\ast_1\om^\ast_2\om^\ast_3\om^\ast_5
+\om^\ast_1\om^\ast_2\om^\ast_4\om^\ast_7
+\om^\ast_1\om^\ast_3\om^\ast_4\om^\ast_6
+\om^\ast_1\om^\ast_5\om^\ast_6\om^\ast_7+\\
& &
\om^\ast_2\om^\ast_3\om^\ast_6\om^\ast_7
+\om^\ast_2\om^\ast_4\om^\ast_5\om^\ast_6
+\om^\ast_3\om^\ast_4\om^\ast_5\om^\ast_7, \\
\end{array} $$

$$ \begin{array}{rcl}
\hat{f}_5 & := &
\om^\ast_1\om^\ast_2\om^\ast_3\om^\ast_4\om^\ast_5\om^\ast_6
+\om^\ast_1\om^\ast_2\om^\ast_3\om^\ast_4\om^\ast_5\om^\ast_7
+\om^\ast_1\om^\ast_2\om^\ast_3\om^\ast_4\om^\ast_6\om^\ast_7+ \\
& &
\om^\ast_1\om^\ast_2\om^\ast_3\om^\ast_5\om^\ast_6\om^\ast_7
+\om^\ast_1\om^\ast_2\om^\ast_4\om^\ast_5\om^\ast_6\om^\ast_7
+\om^\ast_1\om^\ast_3\om^\ast_4\om^\ast_5\om^\ast_6\om^\ast_7+ \\
& &
\om^\ast_2\om^\ast_3\om^\ast_4\om^\ast_5\om^\ast_6\om^\ast_7, \\
\hat{f}_6 & := &
\om^\ast_1\om^\ast_2\om^\ast_3\om^\ast_4\om^\ast_5\om^\ast_6\om^\ast_7 .\\
\end{array} $$

\abs
Actually, it turns out that there is no set of primary invariants 
of $S[W^\ast]^G$ having a strictly smaller degree product,
hence $\hcF$ is optimal in the sense of \cite[Ch.3.3.2]{derksen}.
As $\hat{f}_0\in S[W^\ast]^G_1$, 
using the technique described in Remark \ref{primrem}, we find
an optimal set $\cF=\{f_1,\ld f_6\}\sseq S[W'^\ast]^G$ of 
primary invariants, having degrees $\{2,3,3,4,6,7\}$.

\Abs
Next we compute secondary invariants of $S[W'^\ast]^G$.
Using the Hilbert series given in Proposition \ref{hilprop},
we by Proposition \ref{cmprop} and Remark \ref{hilserrem}
conclude that there is a minimal 
set of $18$ secondary invariants, having degrees
$\{0, 4, 5, 5, 6, 7, 8, 9, 9, 10, 10, 11, 12, 13, 14, 14, 15, 19 \}.$
To find such a set of secondary invariants, we first compute the homogeneous 
components $S[W'^\ast]^G_d$, for $d\leq 7$, using linear algebra 
techniques, see \cite[Ch.3.1]{derksen}, implemented in \MAGMA, 
and then consider products of the homogeneous invariants thus found
having appropriate degrees. Thus we successively generate homogeneous 
elements $\cG:=\{g_1,g_2\ld g_{18}\}\in S[W'^\ast]^G$,
repeatedly using the method described in Remark \ref{secrem} 
to ensure that we have 
$\dim_{\F_2}(\lr{\pi(g_j);j\in\{1\ld k\}}_{\F_2})=k$, 
for $k\in\{1\ld 18\}$.

\abs
To apply the method described in Remark \ref{secrem}, we again consider
the invariant ring $S[W'^\ast]^D$, where $D\leq G$ be a $2$-Sylow 
subgroup of $G$, see the proof of Proposition \ref{cmprop}.
As $S[W'^\ast]^D$ is Cohen-Macaulay, using linear algebra
techniques, implemented in \MAGMA, we obtain the finite presentation
$S[W'^\ast]^D\cong\lr{F_1\ld F_6,G_1\ld G_3|R_1\ld R_3}$ as commutative 
$\F_2$-algebras, where the relations are given as
$$
\begin{array}{rcl}
 R_1 & := &
 (F_1+F_2)^2(F_1F_2F_4+F_3F_5+F_4^2+F_4F_5) + (F_3+F_4)F_5^2 + \\
&& (F_1^3+F_1F_2^2+ F_1F_5+ F_2F_5)\cdot G_1
  +(F_1^2F_2+ F_2^3)\cdot G_2 
  +G_1^2, \\
 R_2 & := &
   (F_1F_2+F_2^2+F_3)F_3F_4 +
   (F_1F_2^2+F_1F_3+F_2^3+F_2F_5)\cdot G_1 + \\
&& (F_1F_2^2+F_1F_5+F_2^3+F_2F_5)\cdot G_2 
  +G_1G_2 + G_3, \\
 R_3 & := & (F_2^2F_6 + F_3^2F_5) + F_2F_3\cdot G_2 + G_2^2, \\
\end{array}$$
where the isomorphism from the finitely presented algebra to 
$S[W'^\ast]^D$ is given by $F_i\mt f'_i$ and $G_j\mt g'_j$, where
$\{f'_1\ld f'_6\}\sseq S[W'^\ast]^D$ and 
$\{g'_1\ld g'_3\}\sseq S[W'^\ast]^D$ 
are as in the proof of Proposition \ref{cmprop}.

\abs
Decomposing the set of primary invariants 
$\cF\sseq S[W'^\ast]^G\sseq S[W'^\ast]^D$ into the $\F_2$-algebra 
generators $\{f'_1\ld f'_6\}\cup\{g'_1\ld g'_3\}$ 
of $S[W'^\ast]^D$, again using linear algebra techniques, 
implemented in \MAGMA,
finally yields the finite presentation
$$ S[W'^\ast]^D/(\sum_{i=1}^6 f_iS[W'^\ast]^D)\cong 
 \lr{F_1\ld F_6,G_1\ld G_3|R_1,\ld R_3,R'_1\ld R'_6} $$
as commutative $\F_2$-algebras, where the additional relations are given as
$$
\begin{array}{rcl}
 R'_1 & := & F_2^2 + F_4 + F_5, \\
 R'_2 & := & F_1F_4 + F_2F_5 + G_1 + G_2, \\
 R'_3 & := & F_2F_4 + G_1, \\
 R'_4 & := & F_1^2(F_1+F_2)^2 \!+\! F_1(F_1+F_2)(F_4+F_5) \!+\! F_3(F_3+F_4) 
             \!+\! F_6 \!+\! F_1\cdot G_1, \\
 R'_5 & := &   F_1^2(F_3^2+F_3F_4+F_6)
             + F_1F_2(F_2^2F_4 F_3F_4+F_4F_5+F_6) + \\
      &    &   F_2^2(F_3F_4+F_3F_5+F_4^2) 
             + (F_3+F_4)(F_5^2+F_6) + F_4^2F_5 + \\
      &    &   (F_1F_3+F_1F_5+F_2^3+F_2F_4)\cdot G_1 + \\
      &    &   (F_1F_2^2+F_1F_3+F_1F_4+F_1F_5)\cdot G_2 
             + G_3, \\
 R'_6 & := & F_1F_6(F_3+F_4) .\\
\end{array}$$

\abs
Finally, we end up with secondary invariants 
$\cG:=\{g_1\ld g_6\}\cup\{g_7\ld g_{18}\}$,
where $g_1=1$ and $\{g_1\ld g_6\}$ have degrees $\{0,4,5,5,6,7\}$, while
$$ \begin{array}{rcrr|}
g_7    & := & g_2^2     & (8),\\
g_8    & := & g_2g_3    & (9),\\
g_9    & := & g_2g_4    & (9),\\
g_{10} & := & g_2g_5    &(10),\\
\end{array}\quad
\begin{array}{rcrr|}
g_{11} & := & g_3g_4    &(10),\\
g_{12} & := & g_2g_6    &(11),\\
g_{13} & := & g_2^3     &(12),\\
g_{14} & := & g_2^2g_3  &(13),\\
\end{array}\quad
\begin{array}{rcrr}
g_{15} & := & g_2^2g_5  &(14),\\
g_{16} & := & g_2g_3g_4 &(14),\\
g_{17} & := & g_2^2g_6  &(15),\\
g_{18} & := & g_2^3g_6  &(19),\\
\end{array}$$
where the bracketed numbers indicate the degrees. Hence in particular 
we conclude that $\{f_1\ld f_6\}\cup\{g_1\ld g_6\}\sseq S[W'^\ast]^G$ 
is a minimal $\F_2$-algebra generating set of $S[W'^\ast]^G$, and 
in particular $S[W'^\ast]^G$ is as an $\F_2$-algebra 
generated by invariants of degree at most $7$.

\Rem
In conclusion we note the following observations.

\bfa
As $W$ is a transitive permutation module of odd degree in 
characteristic $2$, there is a non-zero $G$-invariant quadratic 
form on $W$. Indeed, this form coincides with the primary invariant 
$\hat{f}_1\in S[W^\ast]^G$ of degree $2$. Moreover, we have
$$ f_1=a^\ast e^\ast + b^\ast f^\ast + c^\ast d^\ast 
     + d^\ast e^\ast + d^\ast f^\ast + e^\ast f^\ast 
     + d^{\ast 2} + e^{\ast 2} + f^{\ast 2}
   \in S[W'^\ast]^G .$$
Hence this shows that $D_{W'}$ is an embedding 
$D_{W'}\cn G\ra SO_6^+(\F_2)$, where $SO_6^+(\F_2)$ 
denotes the special orthogonal group of degree $6$ over $\F_2$
of maximal Witt index type, see \cite[p.22]{Atlas}.

\abs\bfb
We have $\{0\}\ra (W'/W'')^\ast\ra W'^\ast\ra W''^\ast\ra\{0\}$,
see Section \ref{gdef}.
Note that this extension does not split: Assume to the 
contrary it splits. Then we also have
$W'\cong W''\oplus W'/W''$, and hence 
$(\F_2)_H^G\cong W\cong W''\oplus W'/W''\oplus\F_2$
is semisimple, and thus $\dim_{\F_2}(\End_G((\F_2)_H^G))=3$. 
By the proof of Proposition \ref{hilprop} we have
$\dim_{\Q_2}(\End_G((\Q_2)_H^G))=2$, a contradiction to 
\cite[Thm.3.11.3]{benson}.

\abs
Still we are tempted to apply the technique described at the 
beginning of Remark \ref{primrem}, yielding an embedding
$S[W'^\ast]^G/((W'/W'')^\ast S[W'^\ast])^G\ra S[W''^\ast]^G$.
The question arises whether this map is still surjective.
While the authors do not see a structural reason why this should 
be the case, from the primary invariants 
$\{f_4,f_5,f_6\}\sseq S[W'^\ast]^G$ 
of degrees $\{4,6,7\}$ we obtain $\{c_0,c_1,c_2\}\sseq S[W''^\ast]^G$,
given as
$$ \begin{array}{rcl}
c_2 & := & 
a^{\ast 4} + b^{\ast 4} + c^{\ast 4}
+ a^{\ast 2}b^{\ast 2} 
+ a^{\ast 2}c^{\ast 2} 
+ b^{\ast 2} c^{\ast 2} + \\
& &
  a^{\ast 2}b^\ast c^\ast 
+ a^\ast b^{\ast 2} c^\ast 
+ a^\ast b^\ast c^{\ast 2}, \\
c_1 & := & 
  a^{\ast 4} b^{\ast 2}
+ a^{\ast 2} b^{\ast 4} 
+ a^{\ast 4} c^{\ast 2} 
+ a^{\ast 2} c^{\ast 4} 
+ b^{\ast 4} c^{\ast 2} 
+ b^{\ast 2} c^{\ast 4} + \\
& & 
  a^{\ast 4} b^\ast c^\ast 
+ a^\ast b^{\ast 4} c^\ast 
+ a^\ast b^\ast c^{\ast 4} 
+ a^{\ast 2} b^{\ast 2} c^{\ast 2}, \\
c_0 & := & 
  a^{\ast 4} b^{\ast 2} c^\ast 
+ a^{\ast 4} b^\ast c^{\ast 2}
+ a^{\ast 2} b^{\ast 4} c^\ast 
+ a^\ast b^{\ast 4} c^{\ast 2} 
+ a^{\ast 2} b^\ast c^{\ast 4} 
+ a^\ast b^{\ast 2} c^{\ast 4}. \\
\end{array} $$
It turns out that these are the {\bf Dickson invariants} of 
$S[W''^\ast]^G$, see \cite[Ch.8.1]{benson2}. 
By Dickson's Theorem, see \cite[Thm.8.1.1]{benson2},
the set $\{c_0,c_1,c_2\}\sseq S[W''^\ast]^G$ is algebraically independent, 
and we have $S[W''^\ast]^G=\F_2[c_0,c_1,c_2]$. 
Hence the above map
$S[W'^\ast]^G/((W'/W'')^\ast S[W'^\ast])^G\ra S[W''^\ast]^G$
indeed is surjective.

\abs\bfc
Using the technique described in Remark \ref{primrem}, we may also 
find an optimal set of primary invariants of the invariant ring 
$S[\wti{W'^\ast}]^G$, which by the Hochster-Eagon Theorem, 
see \cite[Thm.4.3.6]{benson2} or Proposition \ref{cmsubgprop}, 
is Cohen-Macaulay. These optimal primary invariants turn out to have 
degrees $\{2,3,3,4,4,7\}$. Indeed, the Hilbert series 
$H_{S[\wti{W'^\ast}]^G}(t)=H_{S[W'^\ast]^G}(t)\in\Q(t)$,
see the proof of Proposition \ref{hilprop}, can be rewritten as
$$ H_{S[\wti{W'^\ast}]^G}(t)=
   \frac{1+2t^5+2t^6+t^7+t^{10}+2t^{11}+2t^{12}+t^{17}}
        {(1-t^2)\cdot (1-t^3)^2\cdot (1-t^4)^2\cdot (1-t^7)}.$$

\abs
Since the optimal set $\cF$ of primary invariants 
has degrees $\{2,3,3,4,6,7\}$, 
there is no set of primary invariants of $S[W'^\ast]^G$ 
having degrees $\{2,3,3,4,4,7\}$. This shows, although $S[\wti{W'^\ast}]^G$ 
and $S[W'^\ast]^G$ do have the same Hilbert series, that their ring 
structures are different. 
Still, as the underlying modules $W'^\ast$ and $\wti{W'^\ast}$ are 
closely related, the corresponding invariant rings should be closely 
related as well. But how this relationship might look like, for the time
being remains mysterious to the authors.


\abs

\abs
{\sc J. M.: \\
Lehrstuhl D f\"ur Mathematik, RWTH Aachen, \\
Templergraben 64, D-52062 Aachen, Germany}; \\
\url{Juergen.Mueller@math.rwth-aachen.de}.

\abs
{\sc
Institut f\"ur Experimentelle Mathematik (IEM), \\
Universit\"at Duisburg-Essen, \\
Ellernstra{\ss}e 29, D-45326 Essen, Germany}.

\abs
{\sc C. R.: \\
Institut f\"ur Experimentelle Mathematik (IEM), \\
Universit\"at Duisburg-Essen, \\
Ellernstra{\ss}e 29, D-45326 Essen, Germany}.\\
\url{ritzenth@math.jussieu.fr}.

\end{document}